\documentclass[12pt]{article}
\usepackage{amssymb,amsmath}
\usepackage{graphicx}
\usepackage{amsfonts}
\usepackage{float}
\usepackage{flafter}
\usepackage{color}

\labelwidth\leftmargini\advance\labelwidth-\labelsep

\def\R{\relax\ifmmode I\!\!R\else$I\!\!R$\fi}

\def\Z{\relax\ifmmode Z\!\!\!Z\else$Z\!\!\!Z$\fi}

\def\C{\relax\ifmmode C\!\!\!\!I\else$C\!\!\!\!I$\fi}

\def\K{\relax\ifmmode I\!\!K\else$I\!\!K$\fi}

\def\N{\relax\ifmmode I\!\!N\else$I\!\!N$\fi}

\newcounter{defcounter}[section]
{\vspace{0.1cm}\begin{sloppypar}\noindent\stepcounter{defcounter}{\bfseries
Definition
      \thesection.\thedefcounter}}%
{\end{sloppypar}\vspace{0.1cm}}

\newtheorem{corollary}{Corollary}[section]

\newtheorem{theorem}{Theorem}[section]

\newtheorem{proposition}{Proposition}[section]

\newcommand{\proof}{{\bf Proof.} }

\newcommand{\qed}{\hfill $\square$}

\begin{document}
\thispagestyle{empty}
\begin{center}
{\Large {\bf Representations of real numbers induced by probability distributions on $\mathbb{N}$ }}
\end{center}
\begin{center}J\"org Neunh\"auserer\\
Leibnitz University of Hannover, Germany\\
joerg.neunhaeuserer@web.de
\end{center}
\begin{center}
\begin{abstract}
We observe that a probability distribution supported by $\mathbb{N}$, induces a representation of real numbers in $[0,1)$ with digits in $\mathbb{N}$. We first study the Hausdorff dimension of sets with prescribed digits with respect to these representations. Than we determine the prevalent frequency of digits and the Hausdorff dimension of sets with prescribed frequencies of digits. As examples we consider the geometric distribution, the Poisson distribution and the zeta distribution.
\end{abstract}
\end{center}
\section{The Representations}
Let ${\mathfrak p}=(p_{i})_{i\in\mathbb{N}}$ be a probability distribution supported by $\mathbb{N}$, this means $p_{i}\in(0,1)$ for all $i\in\mathbb{N}$ and
\[ \sum_{n=1}^{\infty}p_{n}=1.\]
Let $\widehat{p_{1}}=0$ and
\[ \widehat{p_{n}}=\sum_{i=1}^{n-1}p_{i}\]
if $n\ge 2$. For $n\in\mathbb{N}$ we consider linear contractions $T_{n}:[0,1)\to [0,1)$, given by
\[ T_{n}x=p_{n}x+ \widehat{p_{n}},\]
and introduce the map $\pi_{{\mathfrak p}}:\mathbb{N}^\mathbb{N}\to \mathbb{R}$ by
\[ \pi_{{\mathfrak p}}((n_{j}))=\lim_{j\to\infty}T_{n_{1}}\circ T_{n_{2}}\circ\dots\circ T_{n_{j}}(0)\]
\[  =\widehat{p_{n_{1}}}+\sum_{j=1}^{\infty} p_{n_{1}}\cdots p_{n_{j}}\widehat{p_{n_{j+1}}}.\]
The limit in this expressions exists since the maps $T_{n}$ are contractions, moreover it is easy to see that:
\begin{proposition}
For all probability distributions ${\mathfrak p}$ supported by $\mathbb{N}$, the map $\pi_{\mathfrak{p}}:\mathbb{N}^\mathbb{N}\to [0,1)$ is a bijection.
Moreover the map $\pi_{\mathfrak{p}}$ is continuous, if we endorse $\mathbb{N}^\mathbb{N}$ with the metric
\[ d((n_{i}),(m_{i}))=\sum_{i=1}^{\infty}\delta(n_{i},m_{i})2^{-i},\]
where $\delta(x,y)=0$ if $x=y$ and $\delta(x,y)=1$ otherwise.
   \end{proposition}
\proof If $n_{1}<n_{2}$ we have $p_{n_{1}}+\hat p_{n_{1}}\le p_{n_{2}}$ and hence
\[ T_{n_{1}}([0,1))\cap  T_{n_{2}}([0,1))=[\hat p_{n_{1}},p_{n_{1}}+\hat p_{n_{1}})\cap [\hat p_{n_{2}},p_{n_{2}}+\hat p_{n_{2}})=\emptyset.\]
It follows that $\pi_{\mathfrak{p}}$ is injective. Moreover
\[\bigcup_{n=1}^{\infty}T_{n}([0,1))=[0,p_{1})\cup \bigcup_{n=1}^{\infty}[\sum_{i=1}^{n}p_{i},\sum_{i=1}^{n+1}p_{i})=[0,1)\]
hence $\pi_{\mathfrak{p}}$ is surjective. If  $d((n_{i}),(m_{i}))<2^{-u}$ we have $n_{i}=m_{i}$ for $i=1,\dots ,u$, which implies
\[ |\pi_{\mathfrak{p}}((n_{i})-\pi_{\mathfrak{p}}((m_{i})|<\max\{p_{i}|i\in\mathbb{N}\}^u.\]
This proves that $\pi_{\mathfrak{p}}$ is continuous with respect to $d$.
 \qed~\\~\\
For $x\in[0,1)$ we call the sequence $\pi_{\mathfrak{p}}^{-1}(x)$ in $\mathbb{N}^\mathbb{N}$ the representation of $x$ with respect to the probability distribution ${\mathfrak p}$. The entry of this sequence are the digits of $x$ with respect to the representation, given by ${\mathfrak p}$.~\\~\\
Let us look at three examples. For the geometric distribution on $\mathbb{N}$, given by $p_{i}=(1-p)p^{i-1}$ with $p\in (0,1)$, we obtain
\[ \pi_{{\mathfrak p}}((n_{j}))=(1-p^{n_{1}-1})+\sum_{j=1}^{\infty}(1-p)^jp^{n_{1}+\dots+p_{n_{j}-j}}(1-p^{n_{j+1}-1}).\]
For the Poisson distribution  on $\mathbb{N}$ given by $p_{i}=e^{-\lambda}\lambda^{i-1}/(i-1)!$ with $\lambda>0$, we have
\[ \pi_{{\mathfrak p}}((n_{j}))=e^{-\lambda}\sum_{i=1}^{n_{1}-1}\frac{\lambda^{(i-1)}}{(i-1)!}+\sum_{j=1}^{\infty}e^{-\lambda (j+1)}\frac{p^{n_{1}+\dots+n_{j}-j}}{(n_{1}-1)!\cdots (n_{j}-1)!}\sum_{i=1}^{n_{j+1}-1}\frac{\lambda^{(i-1)}}{(i-1)!}.\]
For the zeta distribution, given by $p_{i}=i^{-s}/\zeta(s)$ on $\mathbb{N}$ with $s>1$, we find
\[ \pi_{{\mathfrak p}}((n_{j}))=\zeta(s)^{-1}\sum_{i=1}^{n_{1}-1}i^{-s}+\sum_{j=1}^{\infty}\zeta(s)^{-(j+1)}(n_{1}\cdots n_{j})^{-s}\sum_{i=1}^{n_{j+1}-1}i^{-s}.\]
As far as we know these representations of real numbers were not considered yet.
\section{Prescribed digits}
Let $D\subseteq \mathbb{N}$ be a set of digits. We are interested in the set $\pi_{{\mathfrak p}}(D^{\mathbb{N}})$ of real numbers which have only digits in $D$ in their representation with respect to a probability distribution ${\mathfrak p}$ on $\mathbb{N}$. It turns out that these sets have Lebesgue measure zero if $D\not=\mathbb{N}$. Thus we study the Hausdorff dimension of these sets. Let us recall that the $d$-dimensional Hausdorff measure of a set $B\subseteq \mathbb{R}$ is given by
\[ \mathfrak{H}^{d}(B)=\lim_{\epsilon \to 0^{+}}\inf\{ \sum_{i=1}^{\infty}(b_{i}-a_{i})^d~|~B\subset \bigcup_{i=1}^{\infty}[a_{i},b_{i}];~\forall i\in\mathbb{N}:(b_{i}-a_{i})\le \epsilon\}\]
and the Hausdorff dimension of $B$ is
\[ \dim_{H}B=\inf\{d~|~H^{d}(B)=0\}=\sup\{d~|~H^{d}(B)=\infty\}.\]
We recommend \cite{[FA1]} or \cite{[PE]} as an introduction to dimension theory. Using the notion of Hausdorff dimension, we obtain:
\begin{theorem}
 Let ${\mathfrak p}=(p_{i})_{i\in\mathbb{N}}$ be a probability distribution supported by $\mathbb{N}$ and $D\subseteq\mathbb{N}$. If $d\ge 0$ is the solution of
 \[ \sum_{i\in D}p_{i}^{d}=1,\]
we have $\dim_{H}\pi_{{\mathfrak p}}(D^{\mathbb{N}})=d$.
\end{theorem}
\proof We have
\[ \bigcup_{i\in D} T_{i}(\pi_{{\mathfrak p}}(D^{\mathbb{N}}))=\bigcup_{i\in D} T_{i}\{\lim_{j\to\infty}T_{n_{1}}\circ T_{n_{2}}\circ\dots\circ T_{n_{j}}(0)~|~n_{j}\in D~\forall j\in\mathbb{N}\}\]
\[=\bigcup_{i\in D} \{\lim_{j\to\infty}T_{i}\circ T_{n_{1}}\circ T_{n_{2}}\circ\dots\circ T_{n_{j}}(0)~|~n_{j}\in D~\forall j\in\mathbb{N}\}=\pi_{{\mathfrak p}}(D^{\mathbb{N}}).\]
This means that $\pi_{{\mathfrak p}}(D^{\mathbb{N}})$ is the attractor of the linear iterated function system $\{T_{i}~|~i\in A\}$ on $[0,1)$, see \cite{[HU]} for finite sets $D$ and \cite{[HF]} for infinite sets. The system fullfills the strong open set condition
$T_{i}((0,1))\cap T_{j}((0,1))=\emptyset$ for $i\not= j$. If A is finite, the result directly follows from the classical work of Moran \cite{[MO]}. If A is infinite it follows from theory of infinite iterated function systems see theorem 3.11 of \cite{[HF]} or \cite{[MU]} for a more general approach.
\qed~\\~\\
As a corollary, we obtain an analogon of Jarnik's \cite{[JA]} classical result on continued fractions.
\begin{corollary}
If $B$ is the set of bounded sequences in $\mathbb{N}^\mathbb{N}$, we have $\dim\pi_{{\mathfrak p}}(B)=1$ for all ${\mathfrak p}$.
\end{corollary}
\proof Since Hausdorff dimension is countable stable
\[ \dim_{H}B=\dim_{H}\bigcup_{k=1}^{\infty}\pi_{{\mathfrak p}}(\{1,\dots,k\}^{\mathbb{N}})=\sup\{\dim_{H}\pi_{{\mathfrak p}}(\{1,\dots,k\}^{\mathbb{N}})|k\in\mathbb{N}\}=1.\]
\qed~\\~\\
In the following we consider $D=\{1,\dots, n\}$ and $\Pi_{{\mathfrak p}}(n)=\pi_{{\mathfrak p}}(\{1,\dots,n\}^{\mathbb{N}})$. If ${\mathfrak p}$ is the geometric distribution with $p\in(0,1)$, we obtain $\dim_{H}\Pi_{{\mathfrak p}}(n)=d$, where $d$ is the solution of
\[ (1-p^{dn})(1-p)^{d}/(1-p^d)=1.\]
We list the first digits of $d$ for some $n$ and $p$ in the following table:
\begin{center}
  \begin{tabular}{ | c || c |c|c|c|c|  }
    \hline
     n/p & 0.1 & 0.25 & 0.5 & 0.75 & 0.9    \\ \hline\hline
    2 & 0.96875 & 0.88920 & 0.69424 & 0.45439 & 0.29434 \\ \hline
    3 & 0.99718 & 0.97718 & 0.87914 & 0.66352 & 0.45656 \\ \hline
    4 & 0.99972 & 0.99463 & 0.94677 & 0.77979 & 0.56428 \\ \hline
    5 & 0.99997 & 0.99868 & 0.97522 & 0.85084 & 0.64218 \\ \hline
    6 & 0.99999 & 0.99967 & 0.98810 & 0.89611 & 0.70137 \\ \hline
  \end{tabular}
\end{center}
For all $n \ge 2$ we have
\[ \lim_{p\to 0}(1-p^{1\cdot n})(1-p)^{1}/(1-p^1)=1\mbox{ and } \lim_{p\to 1} (1-p^{dn})(1-p)^{d}/(1-p^d)=0~\forall d>0,\]
hence the dimension attains all values in $(0,1)$ for $p\in(0,1)$ by continuity.\\ \\
Now let ${\mathfrak p}$ be the Poisson distribution with $\lambda>0$.  We have $\dim_{H}\Pi_{{\mathfrak p}}(n)=d$, where $d$ is the solution of
\[ \sum_{i=1}^{n}e^{-\lambda d}\lambda^{d(i-1)}/((i-1)!)^d=1.\]
We again list the first digits of $d$ for some $n$ and $\lambda$ in a table:
\begin{center}
  \begin{tabular}{ | c || c |c|c|c|c|  }
    \hline
     $n/\lambda$ & 0.25 & 0.5 & 1 & 2 & 4    \\ \hline\hline
    2 & 0.94980 & 0.87189 & 0.69314 & 0.42577 & 0.21288 \\ \hline
    3 & 0.99642 & 0.98345 & 0.92666 & 0.73178 & 0.40665 \\ \hline
    4 & 0.99978 & 0.99809 & 0.98458 & 0.89758 & 0.59553 \\ \hline
    5 & 0.99998 & 0.99981 & 0.99715 & 0.96598 & 0.75770 \\ \hline
    6 & 0.99999 & 0.99998 & 0.99954 & 0.98989 & 0.87203 \\ \hline
  \end{tabular}
\end{center}
For all $n \ge 2$ we have
\[ \lim_{\lambda\to 0} \sum_{i=1}^{n}e^{-\lambda \cdot 1}\lambda^{1\cdot(i-1)}/((i-1)!)^1=1\mbox{ and } \lim_{\lambda\to \infty}  \sum_{i=1}^{n}e^{-\lambda d}\lambda^{d(i-1)}/((i-1)!)^d=0~\forall d>0,\]
hence the dimension attains here all values in $(0,1)$ for $\lambda\in(0,\infty)$ by continuity.~\\~\\
Let ${\mathfrak p}$ now be the zeta distribution with $s>0$.  We have $\dim_{H}\Pi_{{\mathfrak p}}(n)=d$, where $d$ is the solution of
\[ \zeta(s)^{-d}\sum_{i=1}^{n}i^{-sd}=1.\]
We list the first digits of $d$ for some $n$ and $s$:
\begin{center}
  \begin{tabular}{ | c || c |c|c|c|c|  }
    \hline
     $n/s$ & 1.5 & 2 & 3 & 4 & 5    \\ \hline\hline
    2 & 0.48999 & 0.66938 & 0.85250 & 0.92844 & 0.96292 \\ \hline
    3 & 0.64468 & 0.80840 & 0.93681 & 0.97675 & 0.99085 \\ \hline
    4 & 0.72165 & 0.86713 & 0.96462 & 0.98947 & 0.99667 \\ \hline
    5 & 0.76813 & 0.89903 & 0.97731 & 0.99433 & 0.99850 \\ \hline
    6 & 0.79946 & 0.91890 & 0.98418 & 0.99659 & 0.99923 \\ \hline
  \end{tabular}
\end{center}
For all $n \ge 2$ we have
\[ \lim_{s\to 1}\zeta(s)^{-d}\sum_{i=1}^{n}i^{-sd}=0~\forall d\mbox{ and } \lim_{s\to \infty} \zeta(s)^{-1}\sum_{i=1}^{n}i^{-s\cdot 1}=1,\]
hence the dimension attains here all values in $(0,1)$ for $s\in(1,\infty)$ as well.
\section{Frequency of digits}
Let ${\mathfrak f}_{{\mathfrak p}}(x,n)$ be the frequency of the digit $n\in\mathbb{N}$ in the representation of $x\in [0,1)$, given by a probability distribution ${\mathfrak p}$, this means
\[ {\mathfrak f}_{{\mathfrak p}}(x,n)=\lim_{i\to\infty}\mbox{Card}\{j| {\pi_{{\mathfrak p}}^{-1}(x)}_{j}=n~|~j=1,\dots, i\}/i,\]
if the limit exists. As expected we have
\begin{theorem}
Let ${\mathfrak p}$ be a probability distribution supported by $\mathbb{N}$. For almost all $x\in [0,1)$ and all $n\in\mathbb{N}$ we have
${\mathfrak f}_{{\mathfrak p}}(x,n)=p_{n}$.
\end{theorem}
\proof Let $\sigma:[0,1)\to[0,1)$ be the piecewise linear expanding map, given by $T^{-1}_{n}$ on $T_{n}([0,1))$ for $n\in\mathbb{N}$.
The measure ${\mathfrak p}$ on $\mathbb{N}$ induces a Bernoulli measure $b$ on $\mathbb{N}^{\mathbb{N}}$. It is well known (and easy to prove) that this measure is ergodic with respect to the shift map $s(n_{k})=n_{k+1}$ on $\mathbb{N}^{\mathbb{N}}$. We refer here to \cite{[KH]} or \cite{[WA]} for introduction to ergodic theory. The map $\pi_{\mathfrak p}$  projects $b$ to the Lebesgue measure $\ell$ on $[0,1)$, $\ell=b\circ\pi_{\mathfrak p}^{-1}$. Since $b$ is ergodic with respect to $s$ and $\sigma\circ\pi_{\mathfrak p}=\pi_{\mathfrak p}\circ s$, the Lebesgue measure $\ell$ is ergodic with respect to $\sigma$. Applying Birkoff's ergodic theorem to characteristic function $\chi_{n}$ of the interval $T_{n}([0,1))=[\hat p_{n},p_{n}+\hat p_{n})$, we obtain
\[ \lim_{i\to\infty}\frac{1}{i}\sum_{j=1}^{i}\chi_{n}(\sigma^{j}(x))=p_{n} \]
for almost all $x\in[0,1)$ with respect to $\ell$. We have $\chi_{n}(\sigma^{j}(x))=1$ if and only if ${\pi_{{\mathfrak p}}^{-1}(x)}_{j}=n$. Hence ${\mathfrak f}_{{\mathfrak p}}(x,n)=p_{n}$.
\qed~\\~\\
This theorem has the following immediate corollary, which reminds us on the classical theory of continued fractions:
\begin{corollary}
For almost all $x$ the representation of $x$ with respect to ${\mathfrak p}$ is unbounded.
\end{corollary}
Now we consider subset of $[0,1)$ with prescribed frequencies of digits with respect to a representation given by ${\mathfrak p}$. Let
${\mathfrak q}=(q_{i})_{i\in\mathbb{N}}$ be a probability distribution on $\mathbb{N}$, not necessary supported by $\mathbb{N}$, this means $q_{i}\in [0,1]$. We define sets with frequencies of digits given by ${\mathfrak q}$ in the following way
\[ F({\mathfrak p},{\mathfrak q})=\{x\in[0,1)~|~{\mathfrak f}_{{\mathfrak p}}(x,n)=q_{n}~\forall n\in\mathbb{N}\}.\]
Recall that the entropy of ${\mathfrak q}$ is
\[ H({\mathfrak q})=-\sum_{i=1}^{\infty}q_{i}\log(q_{i}).\]
provided that the limit exists. Here we set $q_{i}\log(q_{i})=0$ if $q_{i}=0$. See \cite{[WA]} or \cite{[KH]} for an introduction to entropy theory. Moreover let
\[ E(I_{{\mathfrak p}}({\mathfrak q}))=-\sum_{i=1}^{\infty}q_{i}\log(p_{i})\]
provided that the limit exists. This is the expectation of the information content of ${\mathfrak q}$ with respect to ${\mathfrak p}$. With these notations we have
\begin{theorem}
For all probability distributions ${\mathfrak p}$ and ${\mathfrak q}$ on $\mathbb{N}$, where the first distribution is supported by $\mathbb{N}$, we have
\[ \dim_{H} F({\mathfrak p},{\mathfrak q})=H({\mathfrak q})/ E(I_{{\mathfrak p}}({\mathfrak q})),\]
provided that $H({\mathfrak q})$ and  $E(I_{{\mathfrak p}}({\mathfrak q}))$ exists.
\end{theorem}
\proof Let $b$ be the Bernoulli measure, given by ${\mathfrak q}$ on $\mathbb{N}^\mathbb{N}$. Project this measure to $[0,1)$, using $\pi_{{\mathfrak p}}$, $\mu=b\circ\pi_{{\mathfrak p}}^{-1}$. Note that by the law of large numbers we have $\mu(F({\mathfrak p},{\mathfrak q}))=1$. For $x\in F({\mathfrak p},{\mathfrak q})$ let $I_{n_{1}n_{2}\dots n_{k}}(x)$ be the interval of the form $T_{n_{1}}\circ T_{n_{2}}\circ\dots\circ T_{n_{k}}([0,1))$ that contains $x$. By the definition of $F({\mathfrak p},{\mathfrak q})$ we have
\[ \lim_{k\longmapsto \infty} \frac{1}{k}\log \frac{\mu(I_{n_{1}\dots n_{k}}(x))}{\mbox{Length}(I_{n_{1}\dots n_{k}}(x))^s}\]
\[=\lim_{k\longmapsto\infty}\frac{1}{k}\sum_{i=1}^{k}\log(q_{n_{i}})-s\lim_{k\longmapsto\infty}\frac{1}{k}\sum_{i=1}^{k}\log(p_{n_{i}})\]
\[=- H({\mathfrak q})+sE(I_{{\mathfrak p}}({\mathfrak q})),\]
provided that  $H({\mathfrak q})$ and  $E(I_{{\mathfrak p}}({\mathfrak q}))$ exist. In the last equation we used that for $x\in F({\mathfrak p},{\mathfrak q})$ the frequencies of digits in the ${\mathfrak p}$ representation of $x$ is given by
$\mathfrak{q}$. The above equation implies that for all $x\in  F({\mathfrak p},{\mathfrak q})$
\[ \lim_{k\longmapsto \infty}  \frac{\mu(I_{n_{1}\dots n_{k}}(x))}{\mbox{Length}(I_{n_{1}\dots n_{k}}(x))^s}=\{\begin{array}{cc} 0& s<d\\ \infty & s>d \end{array},\]
where
\[ d=H({\mathfrak q})/ E(I_{{\mathfrak p}}({\mathfrak q})).\]
By the local mass distribution principle, see proposition 4.9 of \cite{[FA1]}, we have $\mathfrak{H}^{s}(F({\mathfrak p},{\mathfrak q}))=\infty$ for $s<d$ and $\mathfrak{H}^{s}(F({\mathfrak p},{\mathfrak q}))=0$ for $s>d$. This implies
$\dim_{H}(F({\mathfrak p},{\mathfrak q}))=d$. \qed~\\~\\
As an example we consider the set of numbers  $F({\mathfrak p},{\mathfrak q})$, which have equidistribution digits from $\{1,\dots,n\}$ in their representation, given by ${\mathfrak p}$. In this case ${\mathfrak q}$ is given by $q_{i}=1/n$ for $i=1,\dots,n$. Hence $H({\mathfrak q})=\log(n)$ and
\[ \dim_{H}F({\mathfrak p},{\mathfrak q})=-n\log(n)/\log(p_{1}\cdots p_{n})).\]
For the geometric distribution ${\mathfrak p}$, this gives
\[  \dim_{H}F({\mathfrak p},{\mathfrak q})=\log(n)/((1-n)\log(p)-\log(1-p),\]
where $p\in(0,1)$. For the Poisson distribution ${\mathfrak p}$ with $\lambda >1$ we have
\[   \dim_{H}F({\mathfrak p},{\mathfrak q})=\log(n)/(\lambda-(n-1)\log(\lambda)/2+\sum_{i=1}^{n-1}(n-i)\log(i))\]
and for the zeta distribution ${\mathfrak p}$ with $s>1$ we obtain
\[ \dim_{H}F({\mathfrak p},{\mathfrak q})=\log(n)/(\log(\zeta(s))+s\sum_{i=1}^{n}\log(i)/n).\]

\end{document}